\documentclass[11pt,a4paper]{article}
\pdfoutput=1
\usepackage{amssymb}  
\usepackage{amsthm}
\usepackage{amsmath}
\usepackage{amsfonts}
\usepackage{moreverb}
\usepackage{graphicx}
\usepackage{float}
\usepackage{enumitem}
\usepackage{caption}
\usepackage{subcaption}
\usepackage{authblk}

\title{Active Fault Tolerant Flight Control System Design - A UAV Case Study}
\author[1]{Rudaba Khan}
\author[2]{Paul Williams}
\author[2]{Paul Riseborough}
\author[1]{Asha Rao}
\author[1]{Robin Hill}
\affil[1]{Department of Mathematics and Geospatial Science, RMIT University, Melbourne, Australia.}
\affil[2]{BAE SYSTEMS Australia, Melbourne Australia.}

\begin{document}
\maketitle
\begin{abstract}
Fault tolerance is achieved through multiply redundant hardware systems in large civil aircraft.  This means of achieving fault tolerance is infeasible for small compact unmanned aerial vehicles.  In this paper we apply a fault tolerant control system which exploits analytical redundancy rather than hardware redundancy to an actual UAV model currently in operation via model-in-the-loop simulation.  The fault tolerant control system comprises a nonlinear model predictive controller integrated with an unscented Kalman filter for fault detection and identification.  The results show that our fault tolerant control system design is able to identify engine failure within seconds of fault occurrence and distribute control authority to the healthy actuators to maintain safe flight.   
\end{abstract}

\maketitle

\section*{Nomenclature}
\noindent\begin{tabular}{@{}lcl@{}}
$x_N$  &=& North direction coordinate [m] \\
$x_D$  &=& Down direction coordinate [m]\\
$V_N$  &=& North velocity [m/s] \\
$V_D$  &=& Down velocity [m/s] \\
$\theta$ &=& Pitch angle \\
$q$    &=& Pitch rate\\
$C_b^n$ &=& Body to navigation frame direction cosine matrix \\
$C_n^b$ &=& Navigation to body frame direction cosine matrix \\
$a_X$   &=& Acceleration in body axis, x-direction [m/s$^2$] \\
$a_Z$ &=& Acceleration in body axis, z-direction [m/s$^2$] \\
$V_T$ &=& True airspeed [m/s] \\
$\rho$ &=& Air density \\
$\bar{q}$ &=& Dynamic pressure\\
$u$ &=& Velocity in body frame, x-direction [m/s]\\
$w$ &=& Velocity in body frame, z-direction [m/s]\\
$\alpha$ &=& Angle of attack\\
$\mathbf{V_{\text{wind}}}$ &=& Wind Velocity \\
$\mathbf{\omega_{\text{wind}}}$ &=& Wind turbulence \\
$a_N$ &=& Acceleration in navigation frame, north-direction [m/s$^2$]\\
$a_D$ &=& Acceleration in navigation frame, down-direction [m/s$^2$]\\
$C_m$ &=& Non-dimensional pitching moment coefficient\\
$CX$ &=& Non-dimensional force coefficient in body axis, x-direction\\
$CZ$ &=& Non-dimensional force coefficient in body axis, z-direction\\
$\bar{c}$ &=& Mean aerodynamic chord [m]\\
$S$ &=& Wing area [m$^2$]\\
$b$ &=& Wing span [m]\\
\end{tabular}

\noindent\begin{tabular}{@{}lcl@{}}
$T$ &=& Thrust [N]\\
$v_s$ &=& Speed of sound [m/s]\\
$\delta_{th}$ &=& Throttle\\
$\delta_{e}$ &=& Elevator deflection\\
$\Delta t$ &=& Time step [sec]\\
\end{tabular} 

\section{Introduction}
Fault tolerance in modern flight control is achieved through the design and implementation of multiply-redundant systems, specifically through the addition of supplementary actuators and sensors brought into action in the event of the failure of a member of the principal set of components.  This improves general system reliability and flight safety, but incurs not only the direct cost of added hardware, but also extra weight penalty and additional system complexity. While achieving increased fault tolerance using multiply-redundant hardware is valid for larger conventional aircraft that can physically accommodate this, the compactness of unmanned aerial vehicles (UAVs) makes it impractical and costly.

An alternative solution is the design and implementation of a re-configurable flight controller that can exploit so-called analytical redundancy, as discussed in \cite{Kale2004} and \cite{vos1996application}. Analytical redundancy arises from the existence of inherent redundancies in the system dynamics.  In this paper we show that it is possible to design an active fault tolerant control system for UAVs by exploiting the inherent redundancies in the system dynamics. This analytical redundancy is a much better option for a UAV than hardware redundancy. Our full active fault tolerant control system applied to an actual UAV currently in operation, comprises a nonlinear model predictive controller (NMPC)  along with an unscented Kalman filter for fault detection and identification is able to identify engine failure within seconds of fault occurrence and distribute control authority to the healthy actuators to maintain safe flight. The NMPC controller as developed in \cite{RKPaper1}, uses pseudospectral numerical techniques along with the UKF filter developed in \cite{RKPaper2}. 

In model predictive control, the focus is on designing a controller where the inputs into the controller design are \textit{what to control}, instead of \textit{how to control}.  This is a subtle, but illuminating, difference, meaning inherent system characteristics such as non-linearities and cross-coupling effects can be exploited by the controller, rather than trying to minimize their influence.  

The key to the design of reconfigurable control systems is exploiting the analytical redundancy of the UAV. In this context, the most promising approach to re-configurable and fault-tolerant control is MPC or variants thereof (see \cite{Kale2004}, \cite{Boskovic2001}, \cite{Gopinathan1998}, \cite{maciejowski1999fault}). Predictive control systems are designed by utilising real-time optimization techniques, where a defined objective (or multi-objective) function is optimized subject to plant operational constraints.  

As per the literature fault tolerant flight control has been mainly utilised within the context of large manned aircraft.  Much of the research on UAVs in this context describes the application of FTC to rotorcraft rather than fixed wing aircraft.  Bateman et. al. \cite{bateman2007actuators} believe that gaining airworthiness approval for UAVs in civil airspace requires an increase in reliability and proposes an active FTC system to deal with control surface failures for a UAV.  The fault tolerant control scheme in \cite{bateman2007actuators} consists of a fault detection and identification (FDI) method based on a signal processing approach and a bank of linear quadratic controllers to handle all faults.  In \cite{bateman2008fault} an FTC strategy for the nonlinear model of a UAV equipped with numerous redundant controls is presented.  Here the authors look at asymmetric actuator failures using a sequential quadratic programming (SQP) algorithm that takes into account non-linearities, aerodynamic and gyroscopic couplings, state and control limitations as a means for FTC.

Other applications of FTC for UAV operations can be found in very early papers on the subject by Copeland and Rattan (1994) \cite{copeland1994fuzzy}, Chen et. al (1998) \cite{chen1998lmi} and Wu et. al (1999) \cite{wu1999qft}.  Copeland and Rattan \cite{copeland1994fuzzy} propose a fuzzy logic supervisor algorithm for a reconfigurable flight control law.  The FTC system described by Chen et. al \cite{chen1998lmi}, uses linear matrix inequality to address wing impairment faults on UAVs.  The authors present a multi-objective approach for establishing a matrix inequality formulation.  Wu et. al \cite{wu1999qft} use quantitative feedback theory as the method of choice for FTC on a remote pilotless aircraft.  

In more recent work by Beainy et. al \cite{beainy2009unmanned} FDI is based on neural networks and a reconfigurable controller based on sliding mode control is designed to compensate for the degradation of the actuation on the occurrence of a fault.  Krueger et. al \cite{kruger2012fault}, on the other hand, present an FTC based on an expanded nonlinear model inversion flight control strategy using sliding mode online learning for neural networks.

In \cite{RKPaper3} we applied our FTC system design to a generic fictional aircraft model.  We simulated an engine failure scenario.  The results demonstrated successful fault detection and identification and the re-allocation of control authority to the healthy actuators by the NMPC controller.  In this paper we apply our design to an actual UAV model currently in operation.  Through simulations results we demonstrate that our FTC system design is an effective means of achieving fault tolerance on a real system. 

Section \ref{section:chap6_UAVdetails} summarises the details of the UAV model, followed by an in depth look at the NMPC controller in section \ref{section:chap6_nmpcDev}.  The active FTC system developed in \cite{RKPaper3} is applied to the UAV model in section \ref{section:chap6_UAVactiveFTC} and its effectiveness is demonstrated through investigation of a number of different scenarios.  Finally a conclusion is presented in section \ref{section:chap6_sof}.

\section{Aircraft Details}\label{section:chap6_UAVdetails}
The model used in this case study is of an actual UAV \cite{PWUAVModel}.  The UAV model is of a twin engine, propeller driven aircraft with dimensions as given in Table \ref{table:chap6_UAVData}.

\begin{table}[H]
\caption{UAV Data}
\label{table:chap6_UAVData}
\begin{center}
\begin{tabular}{ll}
\hline
Wingspan, $b$ & $5.5\,m$\\
\hline
Chord, $c$ & $0.55\,m$\\
\hline
Wing Area, $S$ & $3\,m^2$\\
\hline
Mass, $m$ & $36.8\,\text{kg}$\\
\hline
Propeller Diameter, $D$ & $0.4572\,m$\\
\hline
Stall Speed, $V_{\text{stall}}$ & $12\,m/s$\\
\hline
\end{tabular}
\end{center}
\end{table}

The control inputs used to fly the aircraft include the throttle, aileron, elevator, rudder and flaps.  Only longitudinal motion is considered in this paper hence only elevator and throttle inputs are considered. 

The next section provides details of the development of the NMPC controller for the given UAV model.

\section{NMPC Controller Development}\label{section:chap6_nmpcDev}
The UAV model provided uses experimental data with a series of lookup tables incorporated to find force and moment coefficients.  This model forms the plant model to simulate closed loop performance.  The NMPC prediction model (also referred to as the process model) is based on the UAV model but is an approximation model where mathematical expressions are used to find the force and moment coefficients rather than look up tables.  Before the NMPC controller can be developed, the approximation model must be produced, tested and validated after which the NMPC controller can be designed.

The development of the approximation model for the given aircraft is given in the next subsection.

\subsection{Development of Approximation Model}
\begin{figure}[H]
\hspace{-0.8in}
\includegraphics[scale=0.35]{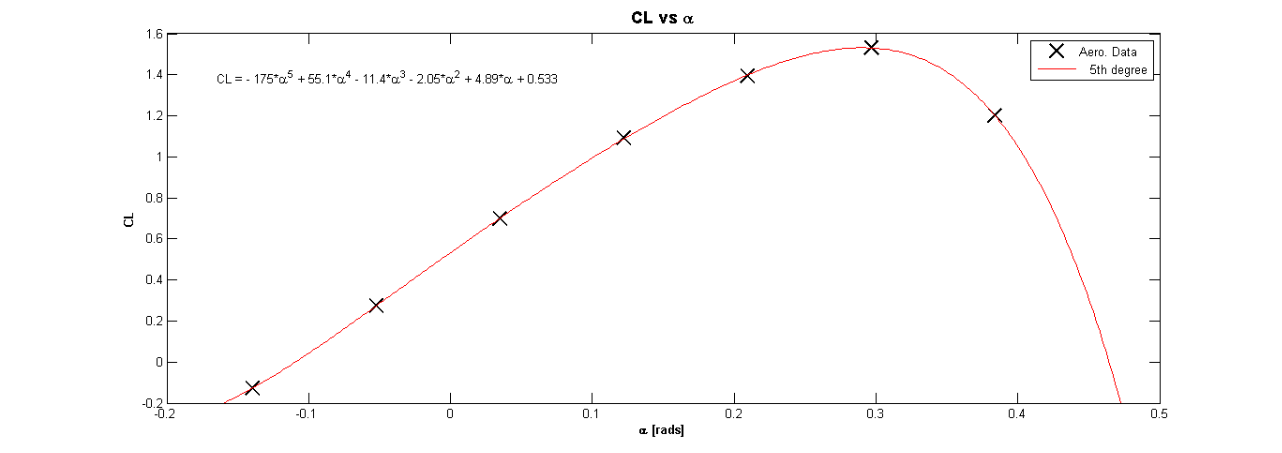} 
%\vspace{-0.5in}
\caption{CL Curve Fitting: experimental data (X), curve fit (red line)}
\label{fig:chap6_CL_alpha}
\end{figure}

\begin{figure}[H]
\hspace{-0.8in}
\includegraphics[scale=0.35]{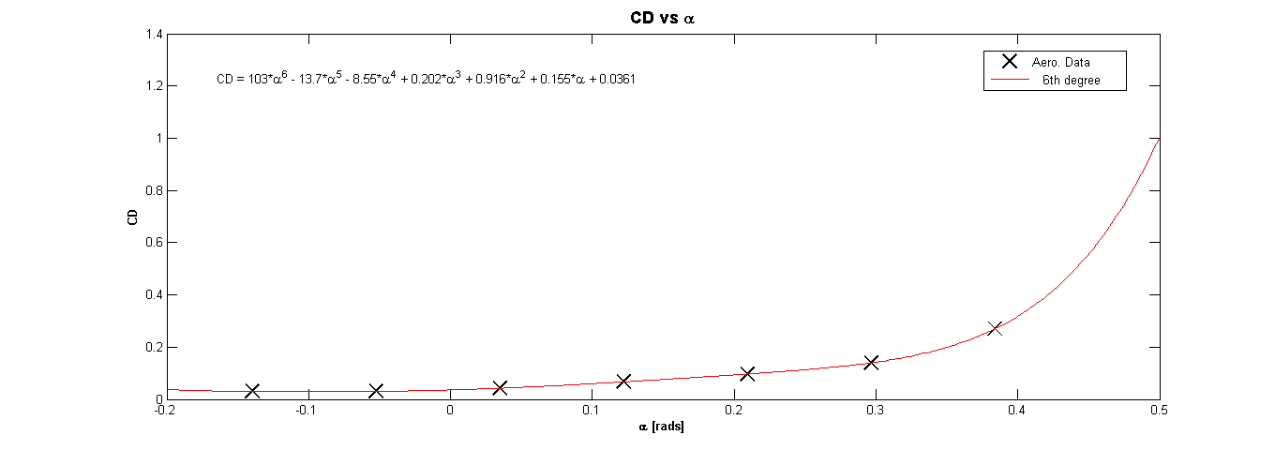} 
%\vspace{-0.5in}
\caption{CD Curve Fitting: experimental data (X), curve fit (red line)}
\label{fig:chap6_CD_alpha}
\end{figure}

\begin{figure}[H]
\hspace{-0.8in}
\includegraphics[scale=0.35]{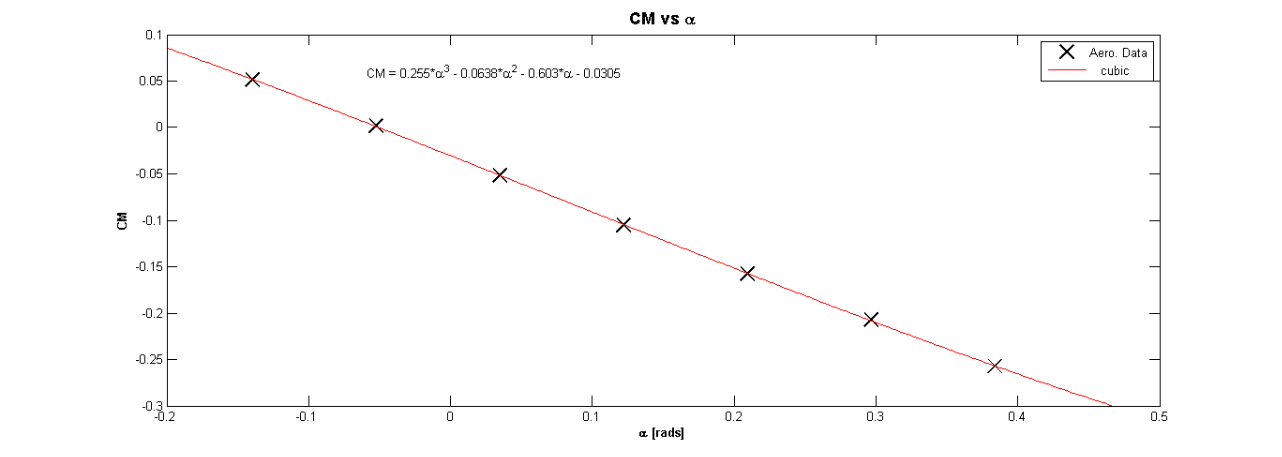} 
%\vspace{-0.5in}
\caption{CM Curve Fitting: experimental data (X), curve fit (red line)}
\label{fig:chap6_CM_alpha}
\end{figure}

The first step in building the approximation model is to fit polynomial curves to the experimental aerodynamic data.  The aerodynamic data includes lift force coefficient $\text{CL}$, drag force coefficient $\text{CD}$ and pitching moment coefficient $\text{CM}$ all of which are given as functions of the angle of attack $\alpha$.  The plots given in figures \ref{fig:chap6_CL_alpha}, \ref{fig:chap6_CD_alpha} and \ref{fig:chap6_CM_alpha} show the curve fits for $\text{CL}$, $\text{CD}$ and $\text{CM}$ respectively.

A fifth order polynomial was fitted to the $\text{CL}$ data using the MATLAB curve fitting toolbox (see figure \ref{fig:chap6_CL_alpha}).  The order of the polynomial is based on trial and error to find the order of best fit.  The following expression was found for $\text{CL}$ as a function of $\alpha$:

\begin{equation}
\text{CL}(\alpha) = - 175\,\alpha^5 + 55.1\,\alpha^4 - 11.4\,\alpha^3 - 2.05\,\alpha^2 + 4.89\,\alpha + 0.533.
\end{equation}

Using the same procedure, sixth order and cubic polynomials were fitted to the $\text{CD}$ and $\text{CM}$ data respectively (figures \ref{fig:chap6_CD_alpha} and \ref{fig:chap6_CM_alpha} respectively).  The following expressions for $\text{CD}$ and $\text{CM}$  were found as functions of $\alpha$:

\begin{equation}
\text{CD}(\alpha) = 103\,\alpha^6 - 13.7\,\alpha^5 - 8.55\,\alpha^4 + 0.202\,\alpha^3 + 0.916\,\alpha^2 + 0.155\,\alpha + 0.0361,
\end{equation}

\begin{equation}
\text{CM}(\alpha) = 0.255\,\alpha^3 - 0.0638\,\alpha^2 - 0.603\alpha \,- 0.0305.
\end{equation}

\begin{figure}[H]
\hspace{-0.8in}
\includegraphics[scale=0.35]{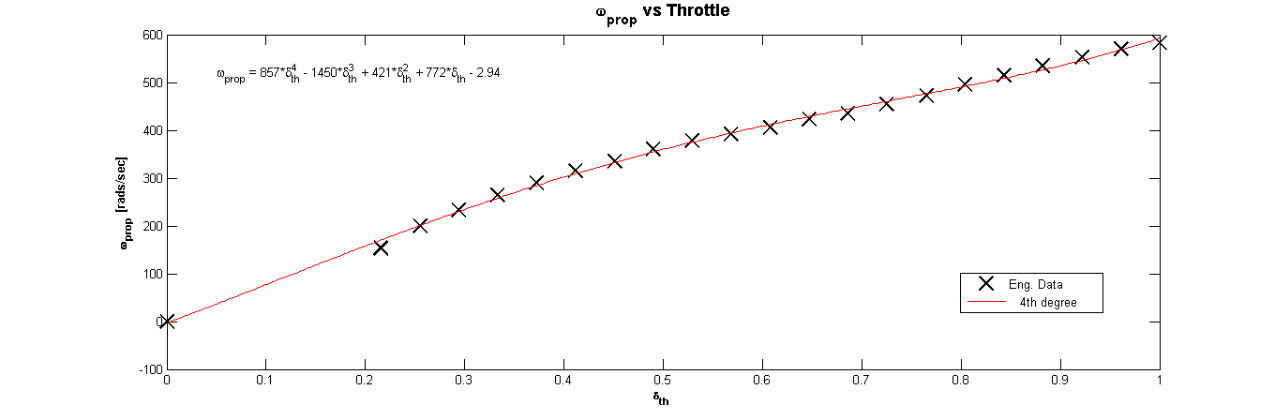} 
%\vspace{-0.5in}
\caption{$\omega_{\text{prop}}$ Curve Fitting: experimental data (X), curve fit (red line)}
\label{fig:chap6_omega_throttle}
\end{figure}

Experimental data was also provided for the engines.  To calculate the thrust force coefficient polynomials were fitted to the angular speed of the propeller ($\omega_{\text{prop}}$ measured in [rads/sec]) vs throttle (figure \ref{fig:chap6_omega_throttle}).

The following expression for $\omega_{\text{prop}}$ as a function of $\delta_{\text{th}}$ was found where $\delta_{\text{th}}$ is the throttle input.

\begin{equation}\label{eqn:chap6_RPM}
\omega_{\text{prop}}(\delta_{\text{th}}) = 857\,\delta_{\text{th}}^4 - 1450\,\delta_{\text{th}}^3 + 421\,\delta_{\text{th}}^2 + 772\,\delta_{\text{th}} - 2.94.
\end{equation}

\begin{figure}[H]
\hspace{-0.8in}
\includegraphics[scale=0.35]{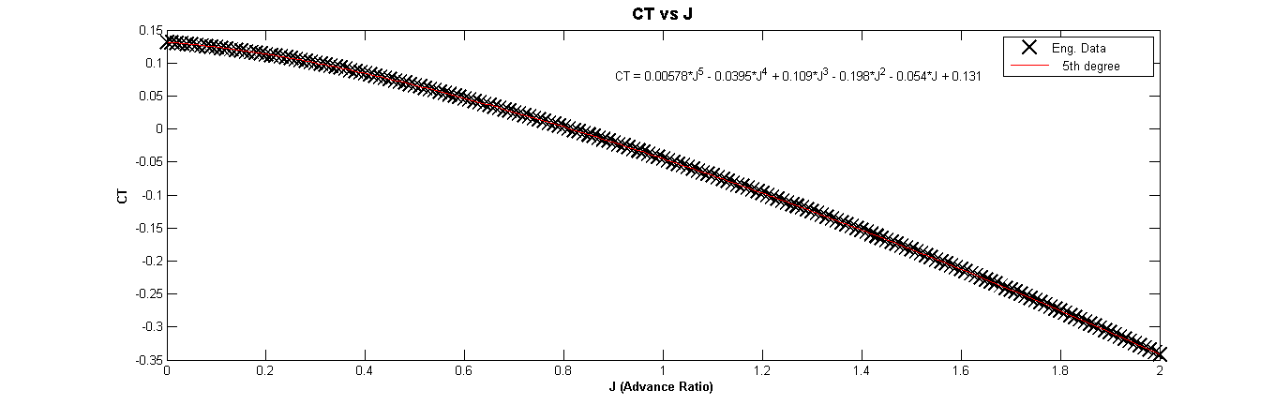} 
%\vspace{-0.5in}
\caption{CT Curve Fitting: experimental data (X), curve fit (red line)}
\label{fig:chap6_CT_J}
\end{figure}

Given an $\omega_{\text{prop}}$ value, the advance ratio, $\mathrm{J}$, of the propeller can be determined from equation \eqref{eqn:chap6_advanceRatio}\cite{mccormick1995}:

\begin{equation}\label{eqn:chap6_advanceRatio}
\mathrm{J}           = \frac{V_T}{n\,D},
\end{equation}

where $n$ is the rotational speed of the propeller given in revolutions per second [RPS]:

\begin{equation}\label{eqn:chap6_prop_n}
n           = \frac{\omega_{\text{prop}}}{2\pi},
\end{equation}

where $\omega_{\text{prop}}$ is calculated by equation \eqref{eqn:chap6_RPM}, $D$ is the propeller diameter and $V_T$ is the true airspeed of the aircraft.  The advance ratio, $J$, is the non-dimensional parameter used to describe the incoming angle of the fluid relative to the propeller blade.  Using experimental data, the thrust coefficient, $CT$ was mapped as a function of $J$, as given in figure \ref{fig:chap6_CT_J}.

A fifth order polynomial was fitted to the data points in figure \ref{fig:chap6_CT_J} and the following expression was found:

\begin{equation}
\text{CT}(\mathrm{J}) = 0.00578\,\mathrm{J}^5 - 0.0395\,\mathrm{J}^4 + 0.109\,\mathrm{J}^3 - 0.198\,\mathrm{J}^2 - 0.054\,\mathrm{J} + 0.131.
\end{equation}

Hence the approximation model of the aircraft for its longitudinal motion is:\\

body velocity components:

\begin{eqnarray}
u &=& V_N\,\cos\theta - V_D\,\sin\theta\,,\\
w &=& V_N\,\sin\theta + V_D\,\cos\theta\,,
\end{eqnarray}

true airspeed:
\begin{equation}
V_T = \sqrt{{V_N}^2 + {V_D}^2}\,,
\end{equation}

angle of attack:

\begin{equation}
\alpha = \arctan\bigg(\frac{w}{u}\bigg)\,,
\end{equation}

dynamic pressure:

\begin{eqnarray}
\bar{q} &=& \frac{1}{2}\,\rho\,\left({V_T}\right)^2,
\end{eqnarray}

aerodynamic forces and moments:
\begin{equation}
\begin{split}
\text{CM} ={}&0.255\,\alpha^3 - 0.0638\,\alpha^2 - 0.603\,\alpha - 0.0305\\
 & + \text{CM}_{\delta_E}\,\delta_E + \text{CM}_q\,\left(\frac{c}{2\,V_T}\right)\,q,
\end{split}
\end{equation}
\\
\begin{equation}
\begin{split}
\text{CL} ={}& - 175\,\alpha^5 + 55.1\,\alpha^4 - 11.4\,\alpha^3 - 2.05\,\alpha^2\\
& + 4.89\,\alpha + 0.533 + \text{CL}_{\delta_E} \,\delta_E + \text{CL}_q\,\left(\frac{c}{2\,V_T}\right)\,q,
\end{split}
\end{equation}
\\
\begin{equation}
\begin{split}
\text{CD} ={}& 103\,\alpha^6 - 13.7\,\alpha^5 - 8.55\,\alpha^4 + 0.202\,\alpha^3\\
& + 0.916\,\alpha^2 + 0.155\,\alpha + 0.0361,
\end{split}
\end{equation}
\\
where $\text{CL}_q = 7.8415$, $\text{CL}_{\delta_E} = 0.3099$, $\text{CM}_q = -14.3808$, $\text{CM}_{\delta_E} = -0.8143$,\\
 
aerodynamic forces and moments in the body axis:
\begin{eqnarray}
\text{CX} &=& \text{CL}\,\sin(\alpha) - \text{CD} \,\cos(\alpha),\\
\text{CZ} &=& -\text{CL}\,\cos(\alpha) - \text{CD}\,\sin(\alpha),\\
\text{FX} &=& \bar{q}\,S\,\text{CX},\\
\text{FZ} &=& \bar{q}\,S\,\text{CZ},\\
\text{M}_p &=& \bar{q}\,S\,c\,\text{CM},\\
\end{eqnarray}

engine forces and moments:
\begin{eqnarray}\label{eqn:chap6_thrustCalcs}
\omega &=& 857\,\delta_{\text{th}}^4 - 1.45e+03\,\delta_{\text{th}}^3 + 421\,\delta_{\text{th}}^2 + 772\,\delta_{\text{th}} - 2.94,\\
\nonumber \\
\text{CT} &=& 0.00578\,J^5 - 0.0395\,J^4 + 0.109\,J^3 - 0.198\,J^2 - 0.054\,J + 0.131,\\
\nonumber \\
\mathbf{F}_{\text{eng}} &=& 2\,\bigg[CT\,\left(\frac{\omega}{2\pi}\right)^2\,\rho\,D^4,\quad 0,\quad 0\bigg]^\intercal,\\
\nonumber \\
\mathbf{M}_{\text{eng}} &=& \mathbf{r}_1 \times \frac{\mathbf{F}_{\text{eng}}}{2} +\mathbf{r}_2 \times \frac{\mathbf{F}_{\text{eng}}}{2},
\end{eqnarray}

where $\mathbf{r}_1$ and $\mathbf{r}_2$ are the position vectors of each engine from the c.g. Note that there are 2 engines hence the multiplication of the engine force by a factor of 2.\\ 

Next we calculate the components of acceleration and the angular rates:
\begin{eqnarray}
a_X &=& \frac{\text{FX}+\text{FX}_{\text{eng}}}{m},\\
\nonumber \\
a_Z &=& \frac{\text{FZ}}{m},\\
\nonumber\\
a_N &=& a_X\,\cos(\theta) + a_Z\,\sin(\theta),\\
\nonumber\\
a_D &=& -a_X\,\sin(\theta) + a_Z\,\cos(\theta) + g.
\end{eqnarray}

Finally the state equations are:
\begin{eqnarray}
\dot{x_D} &=& V_d,\\
\dot{V_N} &=& a_N,\\
\dot{V_D} &=& a_D,\\
\dot{\theta} &=& q\\
\dot{q}   &=& \frac{\text{M}_p+\text{MY}_{\text{eng}}}{I_Y},
\end{eqnarray}

where $I_Y$ is the moment of inertia with respect to the body y-axis and is equal to $18\,\text{kg m}^2$.

In the next section an analysis of the validity of the approximation model is carried out.

\subsection{Model Validation and Verification}
To validate and verify the accuracy of the approximation model a PID controller was developed to control the plant model.  The PID controller consists of a height control loop to calculate a demanded pitch angle given a desired height.  The pitch demand is fed to the pitch control loop to calculate the required elevator deflection and the throttle input is obtained via a speed control loop.

\begin{figure}[H]
\hspace{-0.5in}
\includegraphics[scale=0.8]{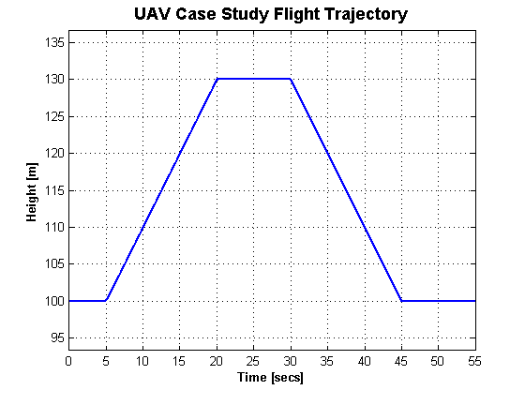} 
\caption{UAV Case Study Reference Trajectory}
\label{fig:chap6_refTraj}
\end{figure}

The PID controller is used to control the plant model which is the full model comprising of the experimental data.  For all validation tests the aircraft is required to fly the trajectory given in figure \ref{fig:chap6_refTraj}.

\begin{figure}[H]
\hspace{-0.8in}
\includegraphics[scale=0.35]{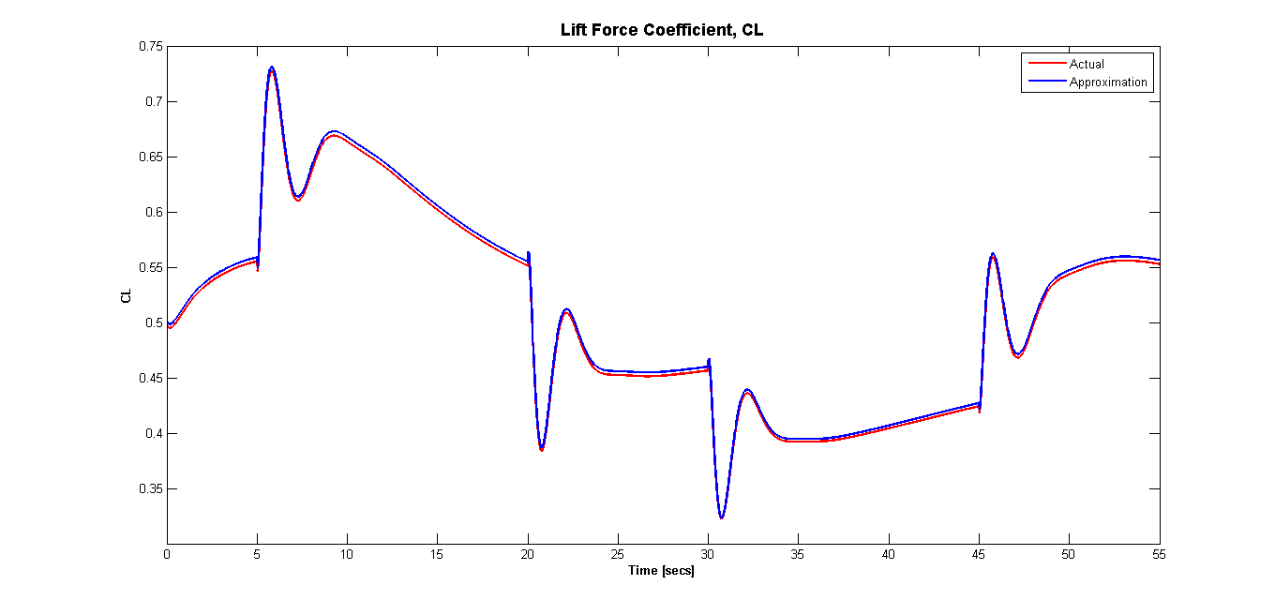} 
%\vspace{-0.5in}
\caption{CL Comparison: actual model (red), approximation model (blue)}
\label{fig:chap6_CL_comparison}
\end{figure}

\begin{figure}[H]
\hspace{-0.8in}
\includegraphics[scale=0.35]{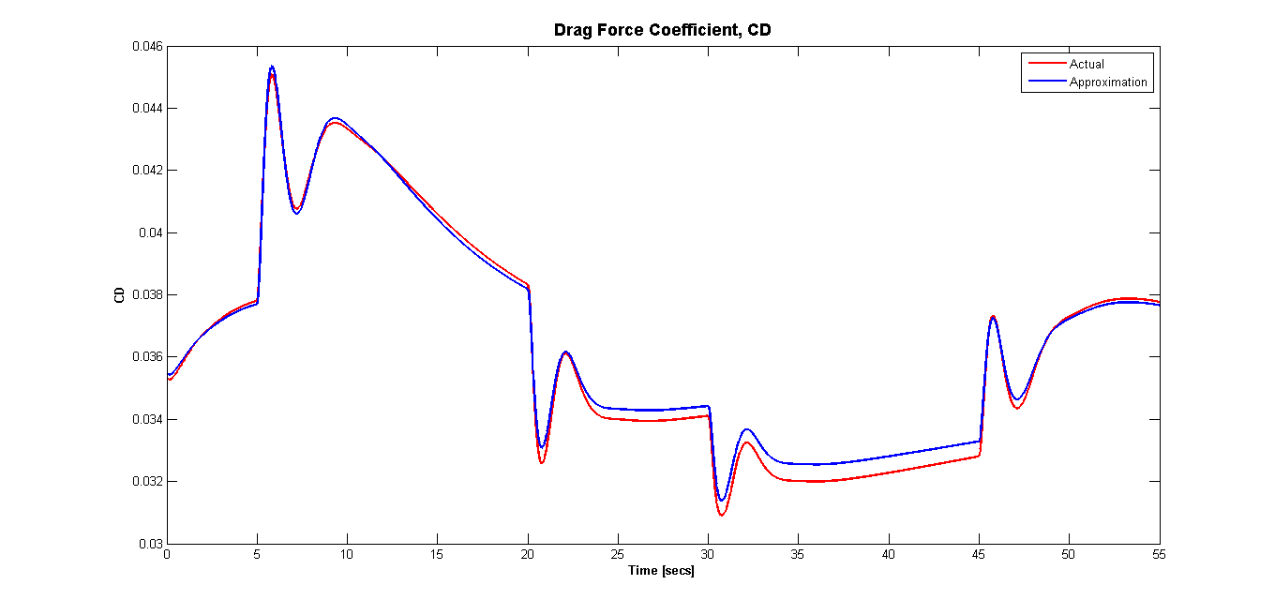} 
%\vspace{-0.5in}
\caption{CD Comparison: actual model (red), approximation model (blue)}
\label{fig:chap6_CD_comparison}
\end{figure}

\begin{figure}[H]
\hspace{-0.8in}
\includegraphics[scale=0.35]{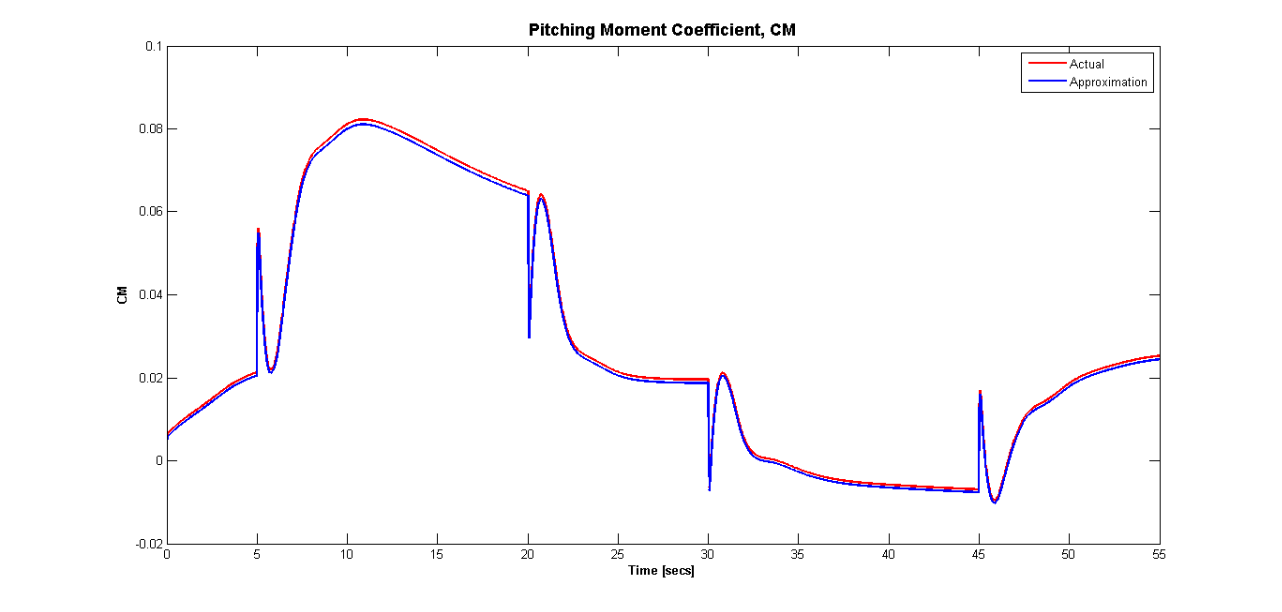} 
%\vspace{-0.5in}
\caption{CM Comparison: actual model (red), approximation model (blue)}
\label{fig:chap6_CM_comparison}
\end{figure}

\begin{figure}[H]
\hspace{-0.8in}
\includegraphics[scale=0.35]{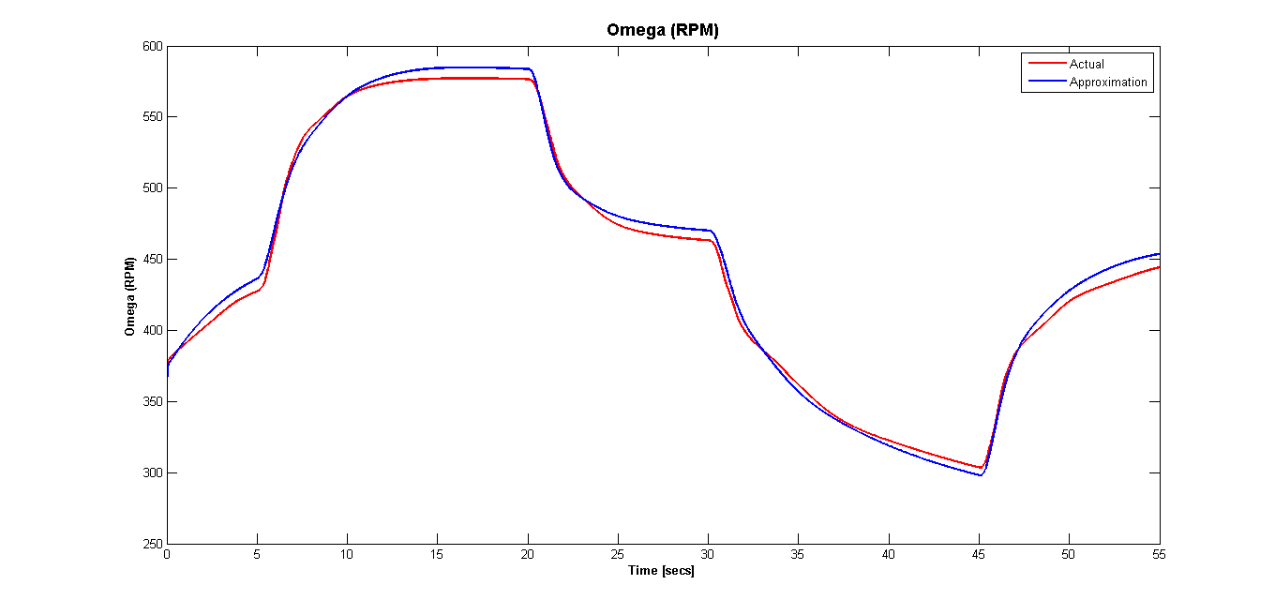} 
%\vspace{-0.5in}
\caption{RPM Comparison: actual model (red), approximation model (blue)}
\label{fig:chap6_Omega_comparison}
\end{figure}

\begin{figure}[H]
\hspace{-0.8in}
\includegraphics[scale=0.35]{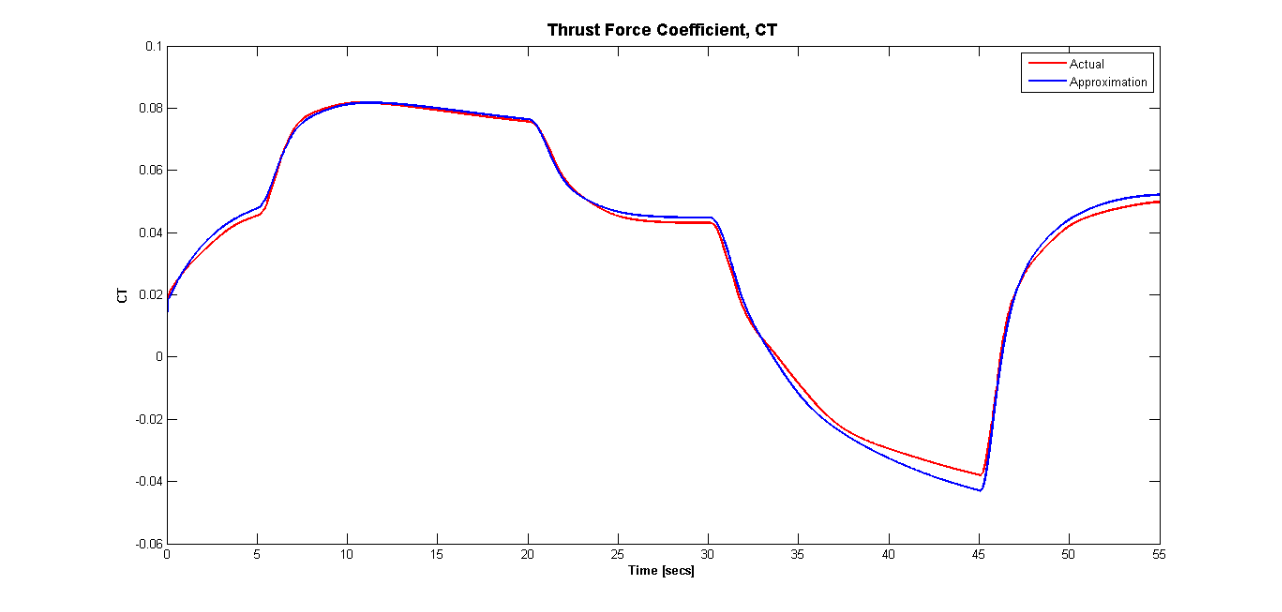} 
%\vspace{-0.5in}
\caption{CT Comparison: actual model (red), approximation model (blue)}
\label{fig:chap6_CT_comparison}
\end{figure}

The first validation was to check if the algebraic expressions for the force and moment coefficients are acceptable.  The control inputs produced by the PID controller to fly the plant model on the given trajectory (figure \ref{fig:chap6_refTraj}) were input into the approximation model and the force and moment coefficients from both models were compared.  Figure \ref{fig:chap6_CL_comparison} shows the $\text{CL}$ values produced by the actual model and the approximation model.  The results show that the approximation model closely follows the response of the actual model.  Figures \ref{fig:chap6_CD_comparison}, \ref{fig:chap6_CM_comparison}, \ref{fig:chap6_Omega_comparison} and \ref{fig:chap6_CT_comparison} show the comparison plots for $\text{CD}$, $\text{CM}$, $\text{RPM}$ and $\text{CT}$ respectively.  All plots show that the approximation model does an excellent job of producing the same results.  There are some differences which are to be expected, however they are within acceptable bounds.  

\begin{figure}[H]
\hspace{-0.8in}
\includegraphics[scale=0.35]{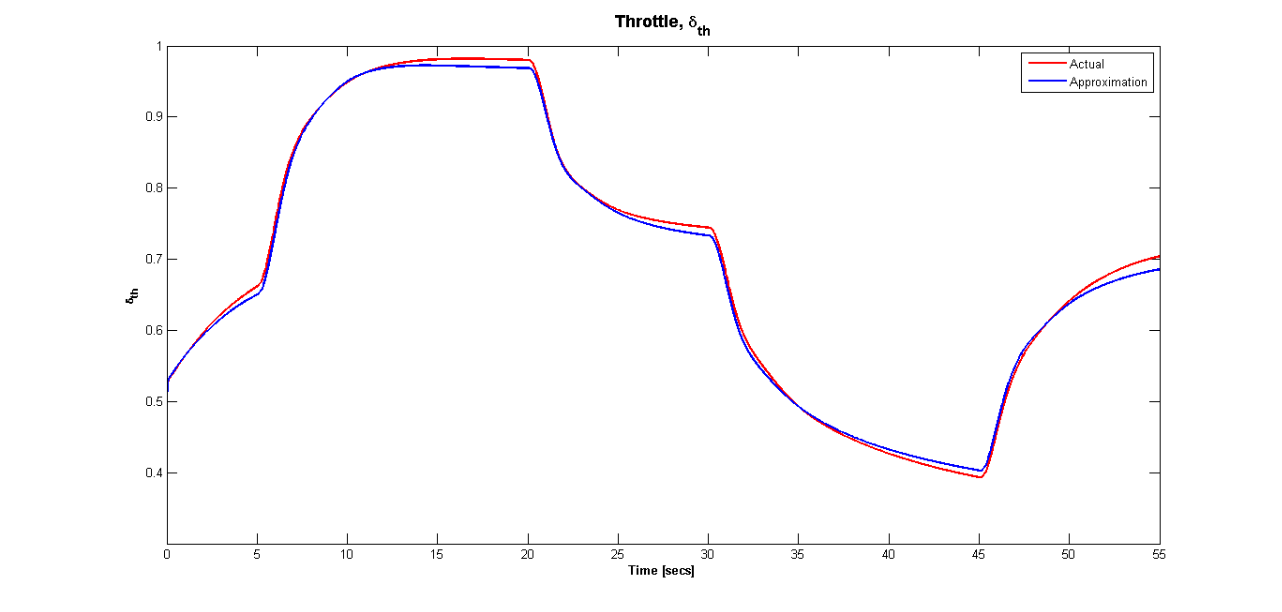} 
%\vspace{-0.5in}
\caption{Throttle Comparison: actual model (red), approximation model (blue)}
\label{fig:chap6_throttle_comparison}
\end{figure}

\begin{figure}[H]
\hspace{-0.8in}
\includegraphics[scale=0.35]{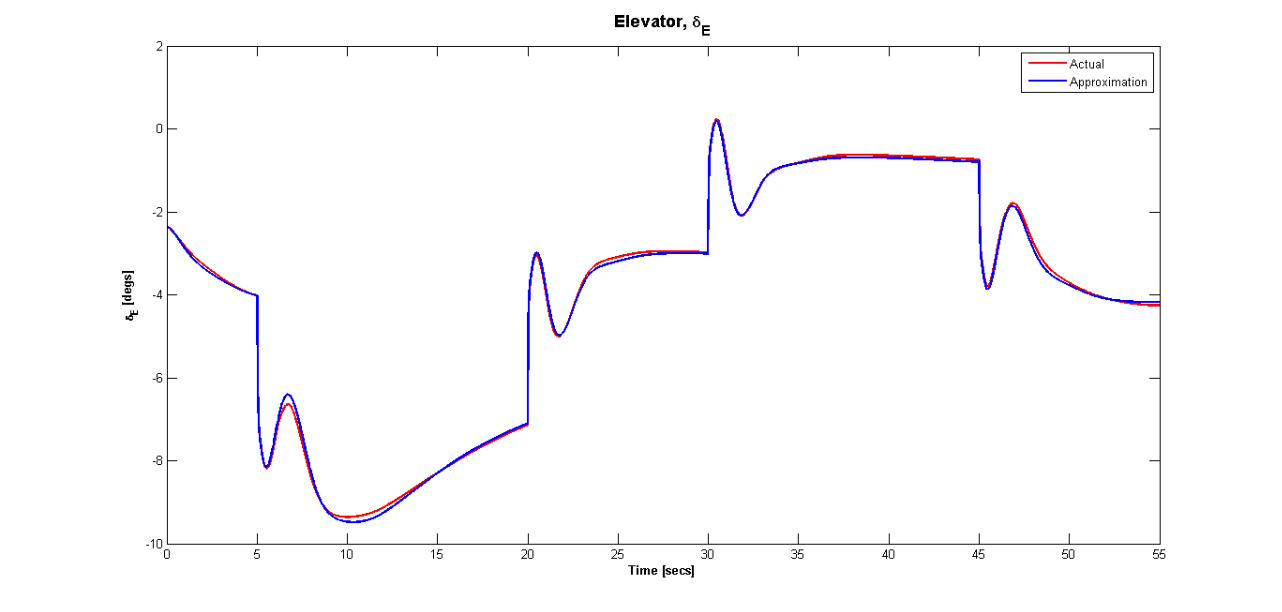} 
%\vspace{-0.5in}
\caption{Elevator Comparison: actual model (red), approximation model (blue)}
\label{fig:chap6_dE_comparison}
\end{figure}

The last step in verifying the validity of the approximation model is to confirm that if the approximation model were to be the plant model then similar control inputs are produced.  The plots given in figures \ref{fig:chap6_throttle_comparison} and \ref{fig:chap6_dE_comparison} show the comparisons of the throttle and elevator responses produced between the actual model and the approximation model.  Both figures show that the approximation is in compliance with the actual model.  

The results produced in this section verify the validity of the approximation model.  This means that the approximation can be used to predict the behaviour of the actual plant within the NMPC controller.  The next section looks at the design of the NMPC controller.

\subsection{NMPC Controller}
A pseudospectral based NMPC controller is designed for the longitudinal motion of the UAV model using the current control inputs, $\delta{th}$ and $\delta_E$.  The state vector $\mathbf{x}$ is defined by:

\begin{equation}\label{eqn:xnmpc}
\mathbf{x} = \left[x_D,\,\,V_N,\,\,V_D,\,\,\theta,\,\,q,\,\,\delta_{\text{th}},\,\,\delta_e,\,\,\Delta{\delta}_\text{th}\,\,\Delta{\delta}_e\right]^\intercal,
\end{equation}

The following optimal control problem is solved at each time step:

\begin{equation}
\begin{split}
\min_{\mathbf{x},\mathbf{u}}\, \frac{H_p}{2}\;\sum_{j = 1}^{j = N+1} & \bigg(\big\Vert \bold{x}_D(j) - \bold{x}_{D_\text{ref}}(j)\big\Vert_{Q_x}^2 + \big\Vert \bold{V}_t(j) - \bold{V}_{t_\text{ref}}(j)\big\Vert_{Q_{VT}}^2 + \big\Vert \bold{V}_D(j) - \bold{V}_{D_\text{ref}}(j)\big\Vert_{Q_{VD}}^2 \\
&\big\Vert \Delta{\delta}_{\text{th}}\big\Vert_{Q_{\text{th}}}^2  + \big\Vert \Delta{\delta}_e\big\Vert_{Q_{\delta_e}}^2 + \big\Vert q\big\Vert_{Q_{q}}^2 + \big\Vert a_D\big\Vert_{Q_{a}}^2\bigg)\;w(j),
\end{split}
\end{equation}

subject to
\begin{eqnarray}
\left(\frac{t_f-t_0}{2}\right)\mathbf{D}_{j,k}\mathbf{x}_j - \mathbf{\dot{x}}_j &=& 0, \\
\mathbf{x}(j_0) - \mathbf{x}_{\text{dem}}(j_0) &=& 0,\\
\mathbf{x}_{lb}  \leq    \mathbf{x}  \leq  \: \mathbf{x}_{ub},\\
\mathbf{u}_{lb}  \leq    \mathbf{u}  \leq  \: \mathbf{u}_{ub},\\
\Delta\mathbf{\delta}_{e_{\text{lb}}}  \leq    \Delta\mathbf{\delta}_e  \leq  \: \Delta\mathbf{\delta}_{e_{\text{ub}}}, \label{eq:chap6_3DOF_cons}
\end{eqnarray}

where $\mathcal{D}_N$ is a spectral differentiation matrix \cite{RKPaper1}, $N$ refers to the number of discretisation (or coincidence) points, $t_0$ and $t_f$ are the initial and final times of the prediction horizon window and the state vector $\mathbf{x}$ is defined in \eqref{eqn:xnmpc} and $\Delta\mathbf{u}$ are the control input rates.  $V_T$ and $V_{T_\text{ref}}$ are the actual and reference true airspeeds respectively, $V_D$ and $V_{D_\text{ref}}$ are the actual and reference down airspeeds respectively and $a_D$ is vertical acceleration in the navigation frame.  $V_T$ and $V_{T_\text{ref}}$ are the actual and reference true airspeeds respectively.  $Q_x$, $Q_{VT}$, $Q_{VD}$, $Q_{th}$, $Q_{\delta_e}$, $Q_q$ and $Q_a$ are diagonal weighting matrices with the following values along the diagonals 5, 1, 1, 0.1, 0.1, 0.01 and 0.01.  The weight values were found through trial and error.  

\begin{table}[H]
\begin{center}
\caption{Constraints for longitudinal Motion}
\begin{tabular}{ccc}
\hline 
\textbf{Variable} &\textbf{ Upper Constraint} & \textbf{Lower Constraint} \\ 
\hline 
$x_D$ & 300 m & 1 m \\ 
\hline 
$V_N$ & $26\,\text{m/s}$  & $15.6\,\text{m/s}$ \\ 
\hline 
$V_D$ & $3\,\text{m/s}$ & $-3\,\text{m/s}$ \\ 
\hline 
$\theta$ & None & None \\ 
\hline 
$q$ & None & None \\ 
\hline 
$\delta_e$ & $30\deg$ & $-30\deg$ \\ 
\hline 
$\delta_{th}$ & $100\%$ & $0\%$ \\ 
\hline
$\Delta{\delta_e}$ & $60\,\text{deg/s}$ & $-60\,\text{deg/s}$ \\ 
\hline 
\end{tabular} 
\end{center}
\end{table}

A prediction window of 5 seconds along with 50 coincidence points were used.  As initial conditions for the controller the trim conditions of the aircraft at a true airspeed of 20m/s are: $\theta = -0.0040 \text{rads}$, $\delta_{th} = 0.7281 \text{rads}$ and $\delta_{E} = -0.0603 \text{rads}$.  To test the controller the aircraft was required to fly the reference trajectory given in figure \ref{fig:chap6_refTraj}.  For the test cases the effects of wind have not been taken into account.\\

\begin{figure}[H]
%\hspace{-1.2in}
\center
\includegraphics[scale=0.35]{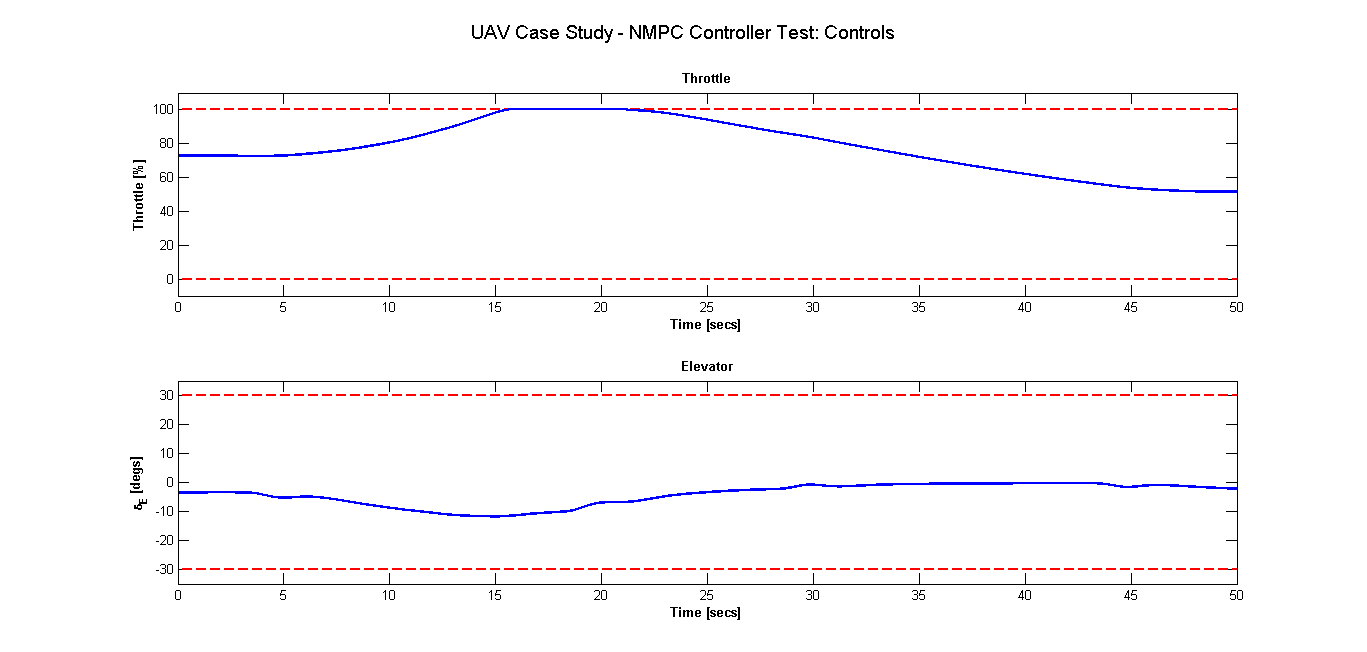} 
%\vspace{-0.5in}
\caption{NMPC Controller Test - Controls, constraints (red), control inputs (blue)}
\label{fig:chap6_UAV_NMPCTest_Controls}
\end{figure}

\begin{figure}[H]
%\hspace{-1.2in}
\center
\includegraphics[scale=0.35]{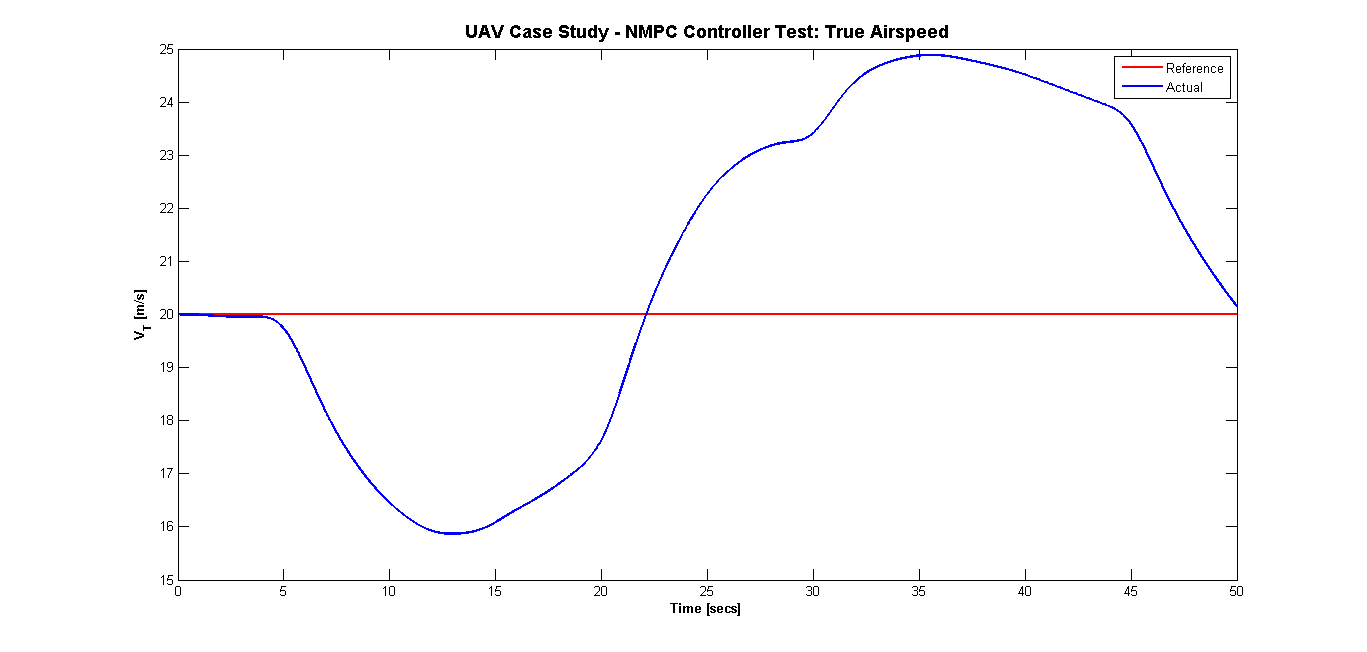} 
%\vspace{-0.5in}
\caption{UAV Case Study NMPC Controller Test - True Airspeed, $V_T$}
\label{fig:chap6_UAV_NMPCTest_VT}
\end{figure}

\begin{figure}[H]
%\hspace{-1.2in}
\center
\includegraphics[scale=0.35]{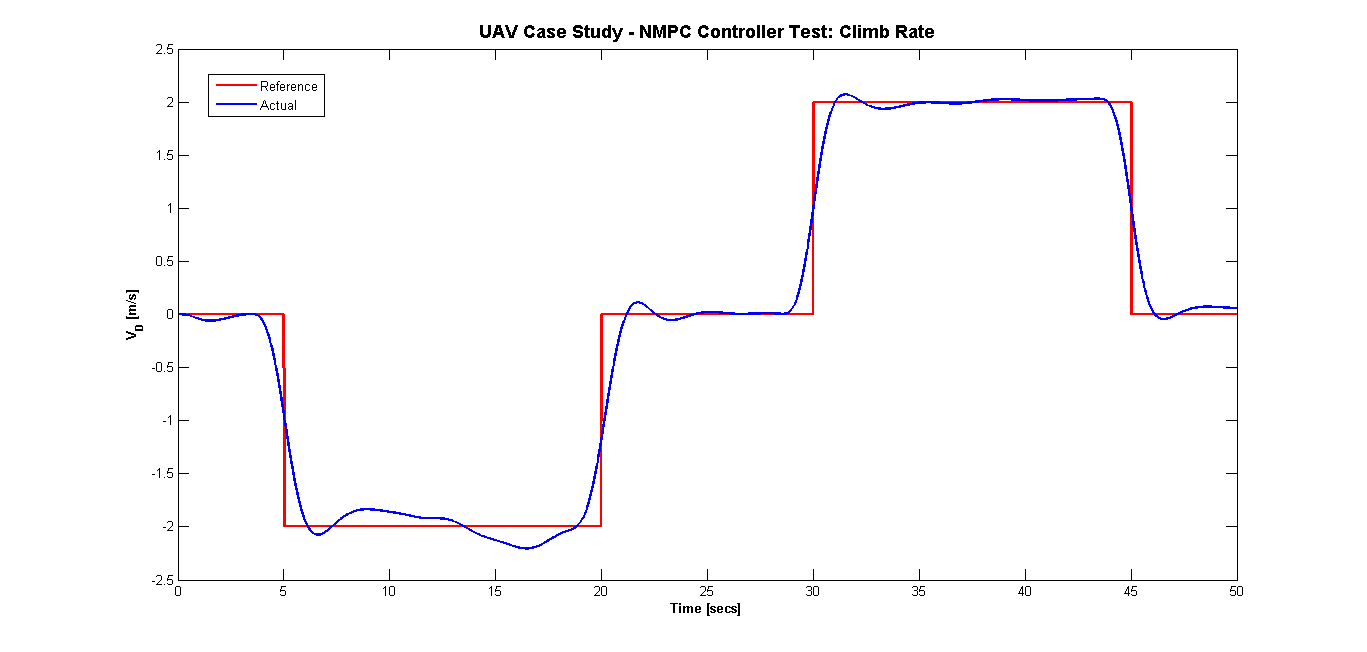} 
%\vspace{-0.5in}
\caption{UAV Case Study NMPC Controller Test - Climb Rate, $V_D$}
\label{fig:chap6_UAV_NMPCTest_VD}
\end{figure}

\begin{figure}[H]
%\hspace{-1.2in}
\center
\includegraphics[scale=0.35]{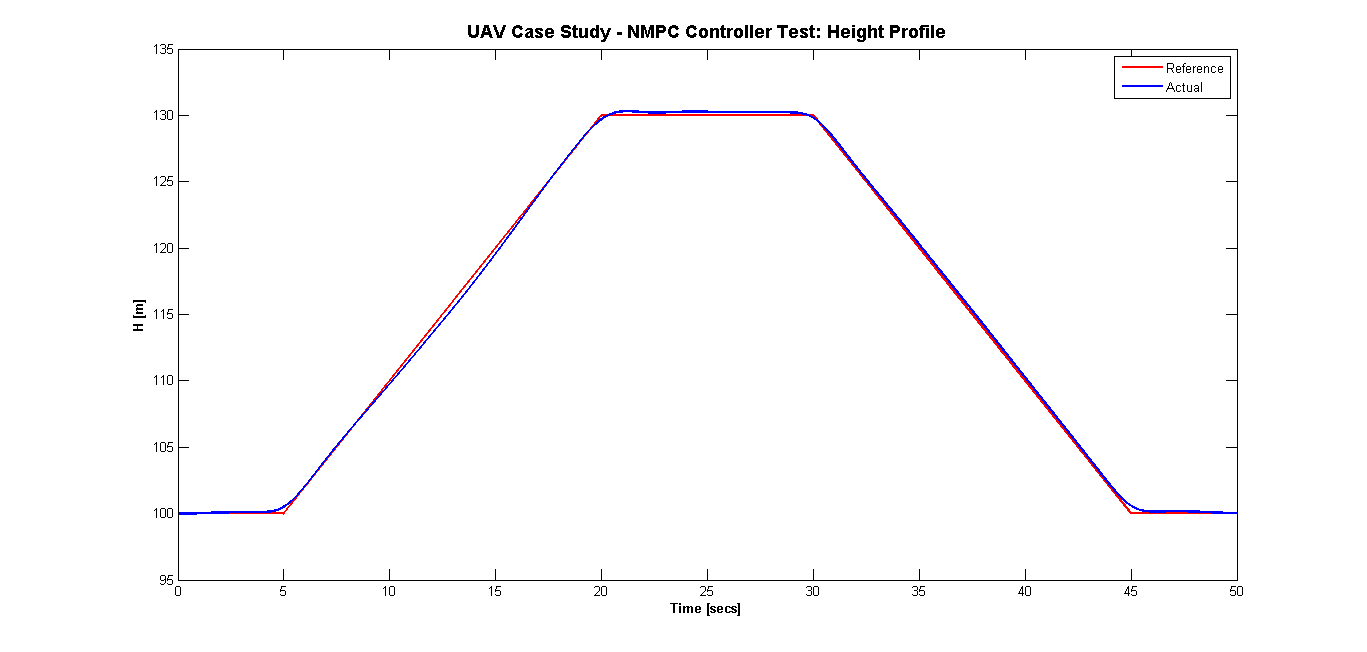} 
%\vspace{-0.5in}
\caption{UAV Case Study NMPC Controller Test - Height Profile}
\label{fig:chap6_UAV_NMPCTest_Height}
\end{figure}

The control inputs produced by the NMPC controller are given in figure \ref{fig:chap6_UAV_NMPCTest_Controls} and show that the inputs remain well within the constraints of the vehicle.  The true airspeed and vertical speed (or climb rate) demands are given in figures \ref{fig:chap6_UAV_NMPCTest_VT} and \ref{fig:chap6_UAV_NMPCTest_VD} respectively.  Both plots show that the controller does an excellent at maintaining the velocity demands.  This is further exemplified by the height profile given in figure \ref{fig:chap6_UAV_NMPCTest_Height} showing the aircraft successfully flying the demanded trajectory.

The next section investigates the application of the active FTC system developed in \cite{RKPaper3} to the given UAV model.

\section{Application of Active FTC System Design to a UAV} \label{section:chap6_UAVactiveFTC}
The active FTC developed in \cite{RKPaper3} is based on a thrust NMPC controller.  Hence the NMPC controller developed in the previous section is now converted to a thrust controller equivalent where the control inputs are now thrust ($\delta_\text{thrust}$) and elevator ($\delta_e$). 

Pseudospectral discretisation is applied to the controller design and the NMPC state vector is:

\begin{equation}\label{eqn:xnmpc_thrust}
\mathbf{x_{nmpc}} = \left[x_D,\,\,V_N,\,\,V_D,\,\,\theta,\,\,q,\,\,\delta_{\text{thrust}},\,\,\delta_e,\,\,\Delta{\delta}_\text{thrust}\,\,\Delta{\delta}_e\right]^\intercal,
\end{equation}

The following optimal control problem is solved each time step:

\begin{equation}
\begin{split}
\min_{\mathbf{x},\mathbf{u}}\, \frac{H_p}{2}\;\sum_{j = 1}^{j = N+1} & \bigg(\big\Vert \bold{x}_D(j) - \bold{x}_{D_\text{ref}}(j)\big\Vert_{Q_x}^2 + \big\Vert \bold{V}_t(j) - \bold{V}_{t_\text{ref}}(j)\big\Vert_{Q_{VT}}^2 + \big\Vert \bold{V}_D(j) - \bold{V}_{D_\text{ref}}(j)\big\Vert_{Q_{VD}}^2 \\
&\big\Vert \Delta{\delta}_{\text{thrust}}\big\Vert_{Q_{\text{T}}}^2  + \big\Vert \Delta{\delta}_e\big\Vert_{Q_{\delta_e}}^2 + \big\Vert q\big\Vert_{Q_{q}}^2 + \big\Vert a_D\big\Vert_{Q_{a}}^2\bigg)\;w(j),
\end{split}
\end{equation}

subject to
\begin{eqnarray}
\left(\frac{t_f-t_0}{2}\right)\mathbf{D}_{j,k}\mathbf{x}_j - \mathbf{\dot{x}}_j &=& 0, \\
\mathbf{x}(j_0) - \mathbf{x}_{\text{dem}}(j_0) &=& 0,\\
\mathbf{x}_{lb}  \leq    \mathbf{x}  \leq  \: \mathbf{x}_{ub},\\
\mathbf{u}_{lb}  \leq    \mathbf{u}  \leq  \: \mathbf{u}_{ub},\\
\Delta\mathbf{u_{\text{lb}}}  \leq    \Delta\mathbf{u}  \leq  \: \Delta\mathbf{u_{\text{ub}}}, 
\label{eq:chap5_6DOF_cons}
\end{eqnarray}

where ${Q_x}$, $Q_{VT}$, $Q_{VD}$, $Q_{T}$, $Q_{\delta_e}$, $Q_{q}$, $Q_{a}$ are weighting matrices with the following values along the diagonal 10, 5, 5, 0.01, 0.1, 0.01 and 0.01 respectively.  The constraints applied are given in Table \ref{table:chap6_constraints_thrust}.

\begin{table}[H]
\begin{center}
\caption{Constraints for UAV Case Study, Thrust Controller}
\label{table:chap6_constraints_thrust}
\begin{tabular}{ccc}
\hline 
\textbf{Variable} & \textbf{Upper Constraint} & \textbf{Lower Constraint} \\ 
\hline 
$x_D$ & 300 m & 1 m \\ 
\hline 
$V_N$ & $26\,\text{m/s}$  & $15.6\,\text{m/s}$ \\ 
\hline 
$V_D$ & $3\,\text{m/s}$ & $-3\,\text{m/s}$ \\ 
\hline 
$\theta$ & None & None \\ 
\hline 
$q$ & None & None \\ 
\hline 
$\delta_e$ & $30\deg$ & $-30\deg$ \\ 
\hline 
$\Delta{\delta}_{\text{thrust}}$ & $122\,\text{N/s}$ & $-122\,\text{N/s}$ \\ 
\hline
$\Delta{\delta}_e$ & $60\,\text{deg/s}$ & $-60\,\text{deg/s}$ \\ 
\hline 
\end{tabular} 
\end{center}
\end{table}

Trim conditions for the aircraft at 20m/s are: $\theta = -0.0040 \text{rads}$, $\text{Thrust} = 27.7426\,N$ and $\delta_{e} = -0.0603 \text{rads}$ and have been used as initial conditions for the controller.  Again a prediction window of 5 seconds along with 50 coincidence points were used.  The upper limit on the thrust constraint continually changes based on the true airspeed of the aircraft according to the equations given in \eqref{eqn:chap6_thrustCalcs}.  The minimum thrust level set always to $0N$.  If however an engine failure is detected the upper limit on thrust is set to the filter estimate plus $2\sigma$ uncertainty.\\ 

A UKF filter was designed to perform FDI with the following process noise and noise covariance matrices:

\begin{equation}
Q = \begin{bmatrix}
(2\Delta t)^2 & 0 & 0 & 0\\
0 & (2\Delta t)^2 & 0 & 0\\
0 & 0 & (0.017 \Delta \, t)^2 & 0\\
0 & 0 & 0 & (122\Delta t)^2
\end{bmatrix},\quad\quad
R = \begin{bmatrix}
(0.5)^2 & 0 & 0\\
0 & (0.5)^2 & 0\\
0 & 0 & (0.17)^2
\end{bmatrix},
\end{equation}

where $\Delta t$ is the filter update rate 0.01 secs.  The initial state vector and covariance matrix are:

\begin{equation}
\mathbf{x}(0) = \left[20,\,\,0\,\,-0.0040,\,\,27.7426\right]^\intercal, \quad
\mathbf{P}(0) = \begin{bmatrix}
(0.5)^2 & 0 & 0 & 0\\
0 & (0.5)^2 & 0 & 0\\
0 & 0 & (0.0850)^2 & 0\\
0 & 0 & 0 & (6)^2
\end{bmatrix}.
\end{equation}

Finally the fault detection logic is based on that given in \cite{RKPaper3}.

\subsection{Numerical Results}
The following scenarios were set up to test the active FTC system on the UAV model:

\begin{enumerate}[label=\bfseries Scenario \arabic*:, leftmargin = 100pt]
\item no fault case
\item engine failure - $50\%$ power loss 20 seconds into flight,
\item engine failure - $70\%$ power loss 30 seconds into flight.
\end{enumerate}

\begin{figure}[H]
\hspace{-0.8in}
\includegraphics[scale=0.35]{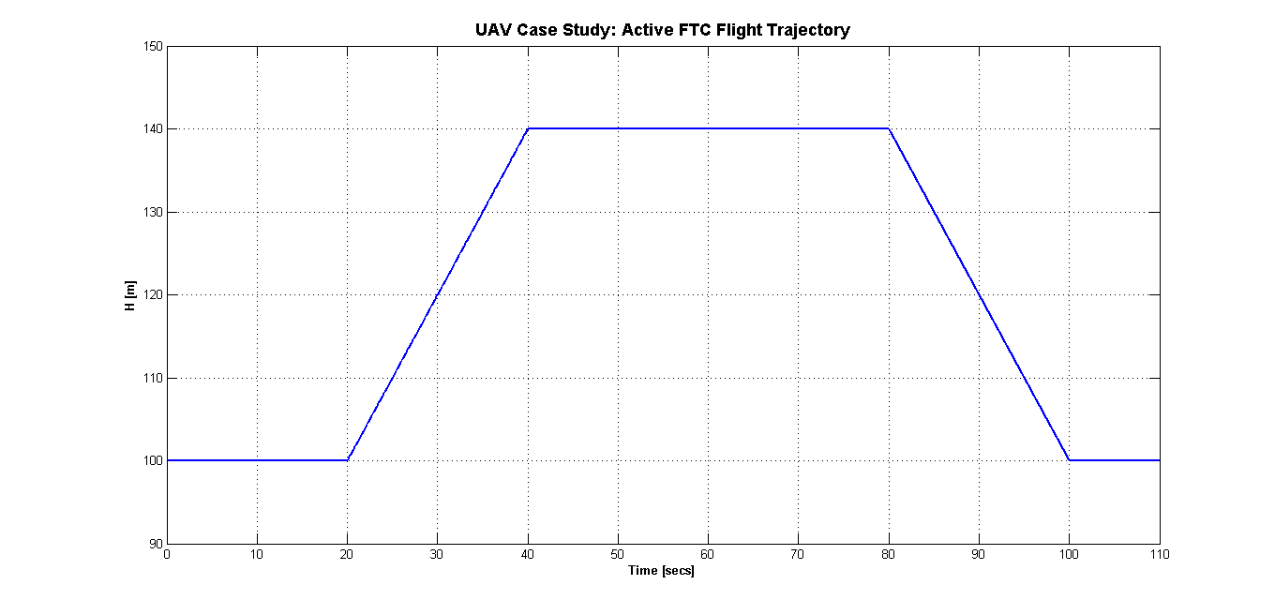} 
%\vspace{-0.5in}
\caption{UAV Case Study Active FTC: Reference Trajectory}
\label{fig:chap6_activeFTC_refTraj}
\end{figure}

The aircraft was required to follow the flight trajectory given in figure \ref{fig:chap6_activeFTC_refTraj} (Note: wind effects have been taken into account using the Dryden Wind Model in MATLAB).

\begin{figure}[H]
\hspace{-0.7in}
\includegraphics[scale=0.35]{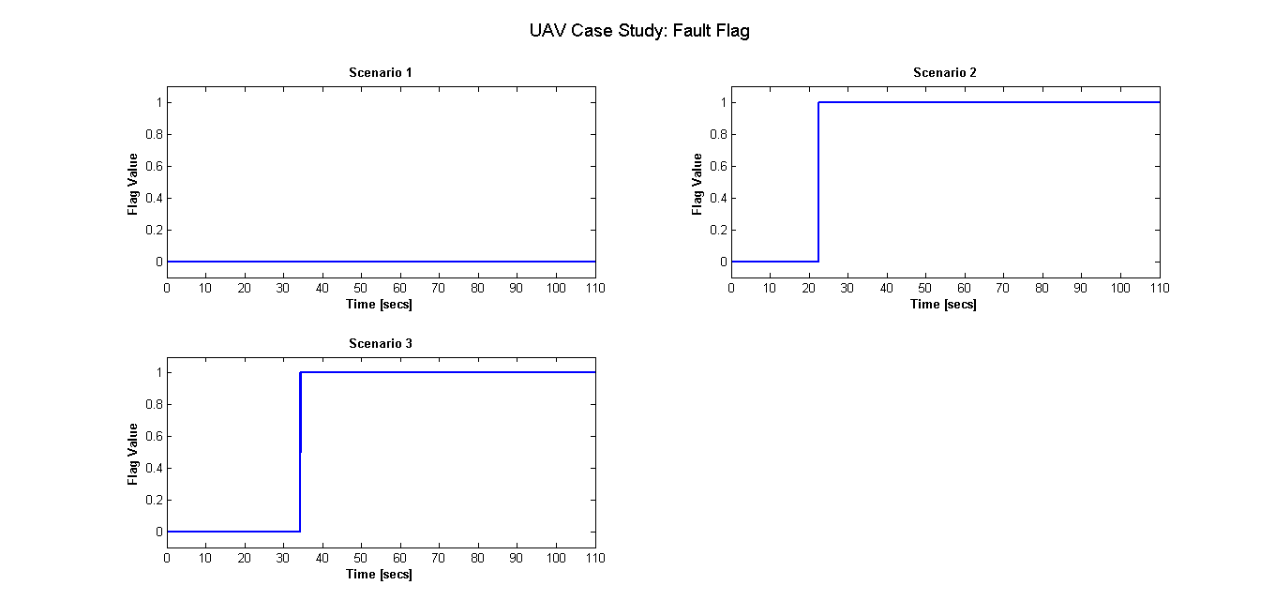} 
%\vspace{-0.5in}
\caption{UAV Case Study Active FTC: Fault Detection Logic Results}
\label{fig:chap6_activeFTC_FaultFlag}
\end{figure}

Figure \ref{fig:chap6_activeFTC_FaultFlag} presents the fault detection logic for each scenario.  There was no fault present in scenario 1 and this is reflected in the results of the fault detection logic as the value of the fault flag remained zero throughout the duration of the flight.  The results for scenarios 2 and 3 show that the fault flag is triggered (ie. the value of the flag switches to 1) within seconds of the fault occurring.

\begin{figure}[H]
%\hspace{-1.2in}
\center
\includegraphics[scale=0.35]{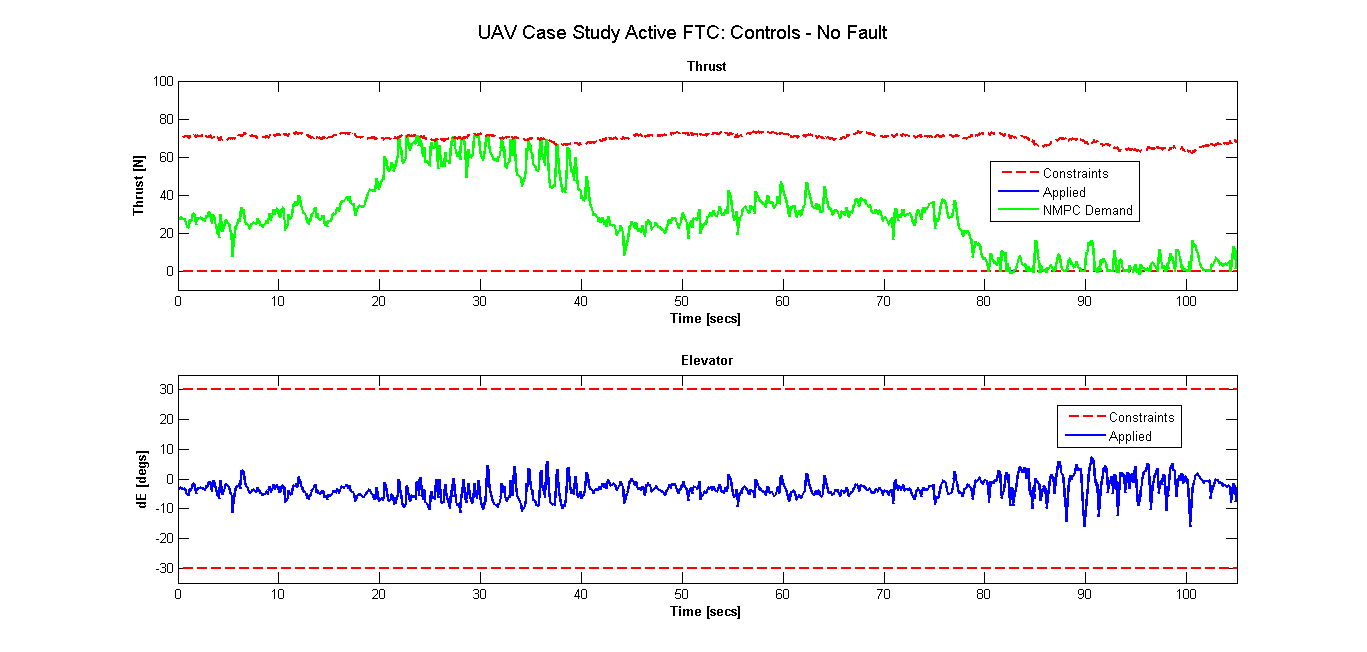} 
%\vspace{-0.5in}
\caption{Scenario 1: No Fault - Controls}
\label{fig:chap6_activeFTC_controls_NF}
\end{figure}

\begin{figure}[H]
%\hspace{-1.2in}
\center
\includegraphics[scale=0.35]{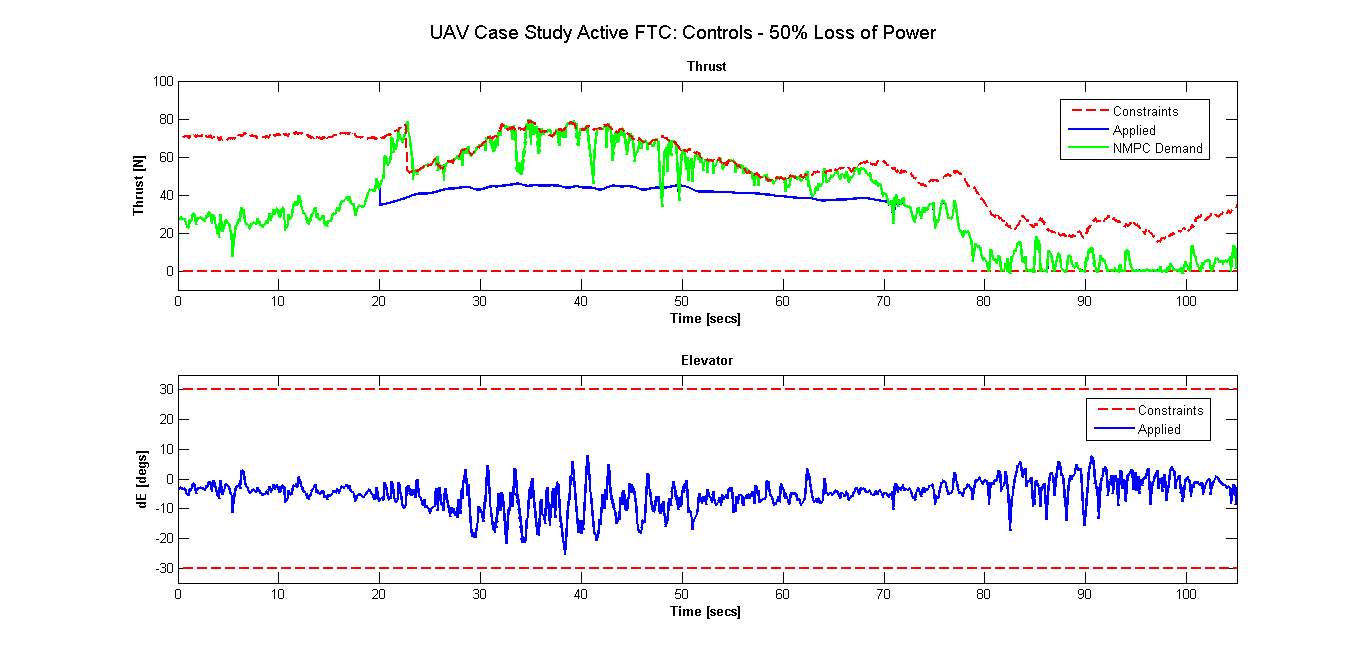} 
%\vspace{-0.5in}
\caption{Scenario 2: $50\%$ Loss of Power - Controls}
\label{fig:chap6_activeFTC_controls_50T}
\end{figure}

\begin{figure}[H]
%\hspace{-1.2in}
\center
\includegraphics[scale=0.35]{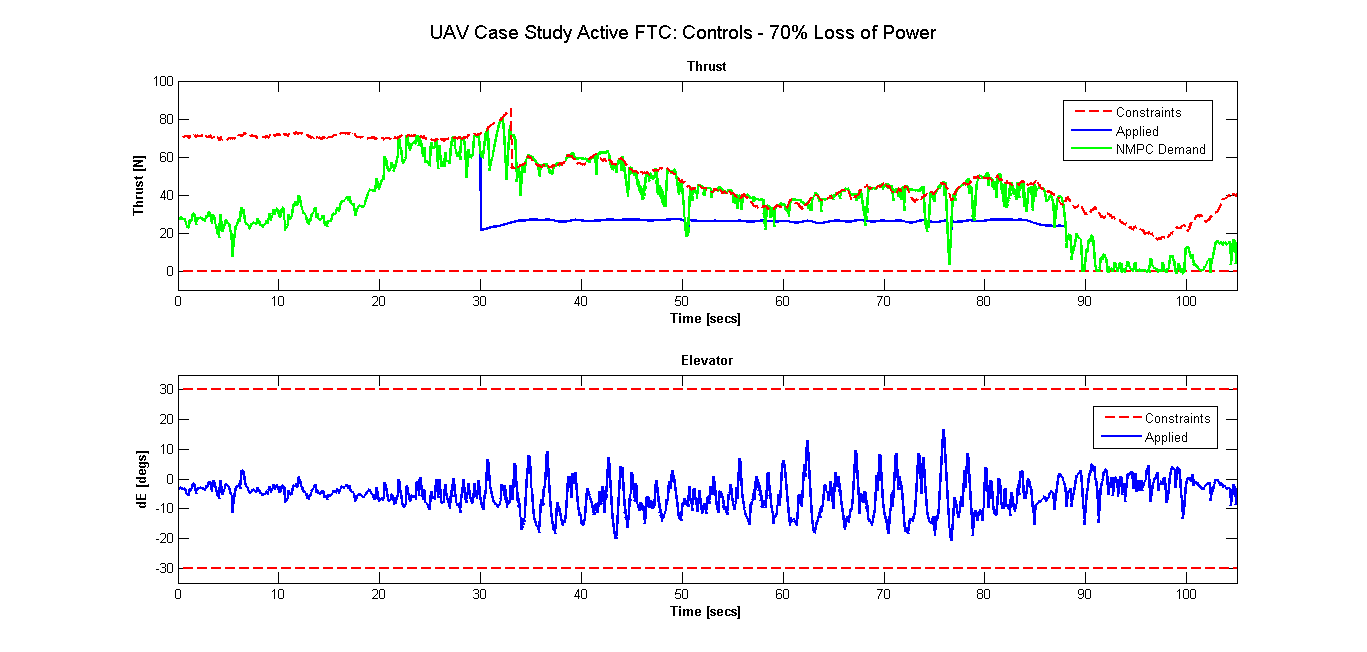} 
%\vspace{-0.5in}
\caption{Scenario 3: $70\%$ Loss of Power - Controls}
\label{fig:chap6_activeFTC_controls_30T}
\end{figure}

The control inputs produced by the controller for each scenario are given in figures \ref{fig:chap6_activeFTC_controls_NF}, \ref{fig:chap6_activeFTC_controls_50T} and \ref{fig:chap6_activeFTC_controls_30T}.  The results clearly indicate that the controller is able to successfully reconfigure based on FDI data and as the loss of power increases the demand on the elevator increases. 

\begin{figure}[H]
%\hspace{-1.2in}
\center
\includegraphics[scale=0.35]{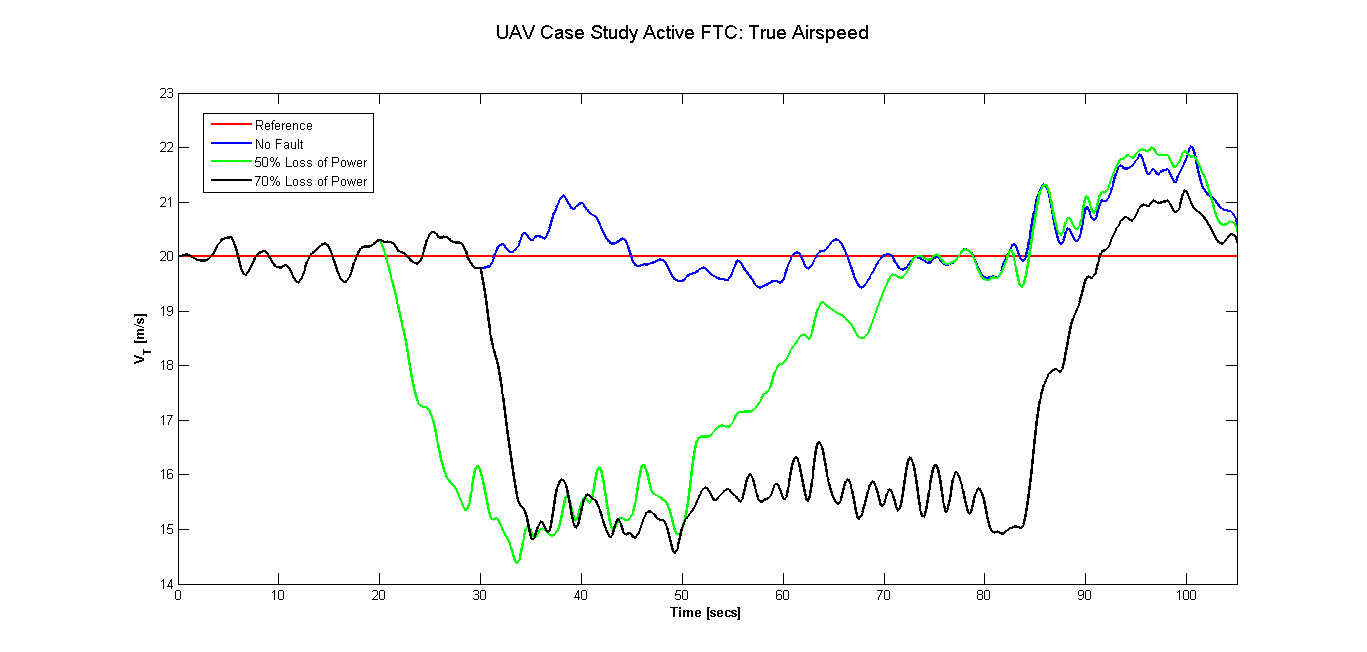} 
%\vspace{-0.5in}
\caption{True Airspeed, $V_T$.  Stall speed did not occur in any scenario.}
\label{fig:chap6_activeFTC_controls_Vt}
\end{figure}

The true airspeed response of the aircraft is given in figure \ref{fig:chap6_activeFTC_controls_Vt}.  The demanded speed was $20m/s$ and the results show that for a $50\%$ loss in power the true airspeed demand is unachievable during climb to altitude.  However during straight and level flight the velocity demand is gradually achieved and is finally reached once the aircraft begins the descent phase.  During the descent the true airspeed demand is similar to the no fault case.  In the case of $70\%$ power loss the aircraft is unable to meet the true airspeed demand during straight and level flight, however halfway through the descent phase the aircraft picks up speed and is able to maintain the reference.  In any given scenario stall speed is never reached.

\begin{figure}[H]
%\hspace{-1.2in}
\center
\includegraphics[scale=0.35]{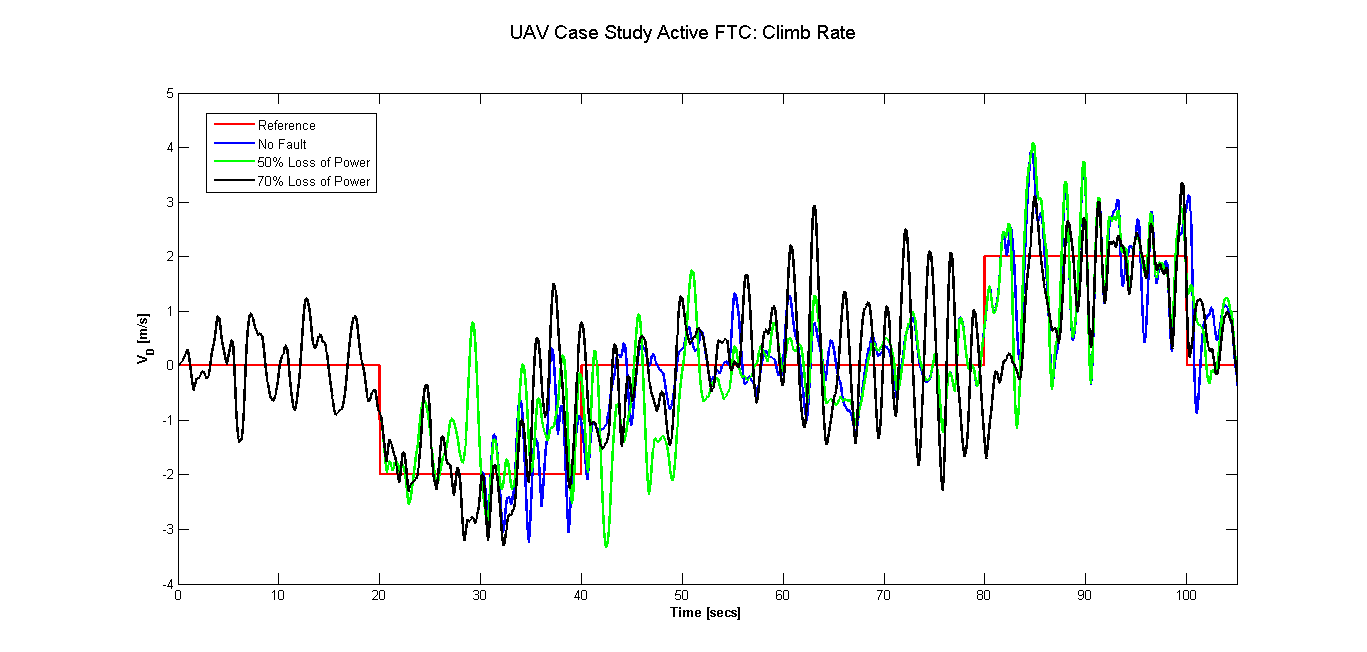} 
%\vspace{-0.5in}
\caption{UAV Case Study Active FTC: Climb Rate, $V_D$}
\label{fig:chap6_activeFTC_controls_Vd}
\end{figure}

\begin{figure}[H]
%\hspace{-1.2in}
\center
\includegraphics[scale=0.35]{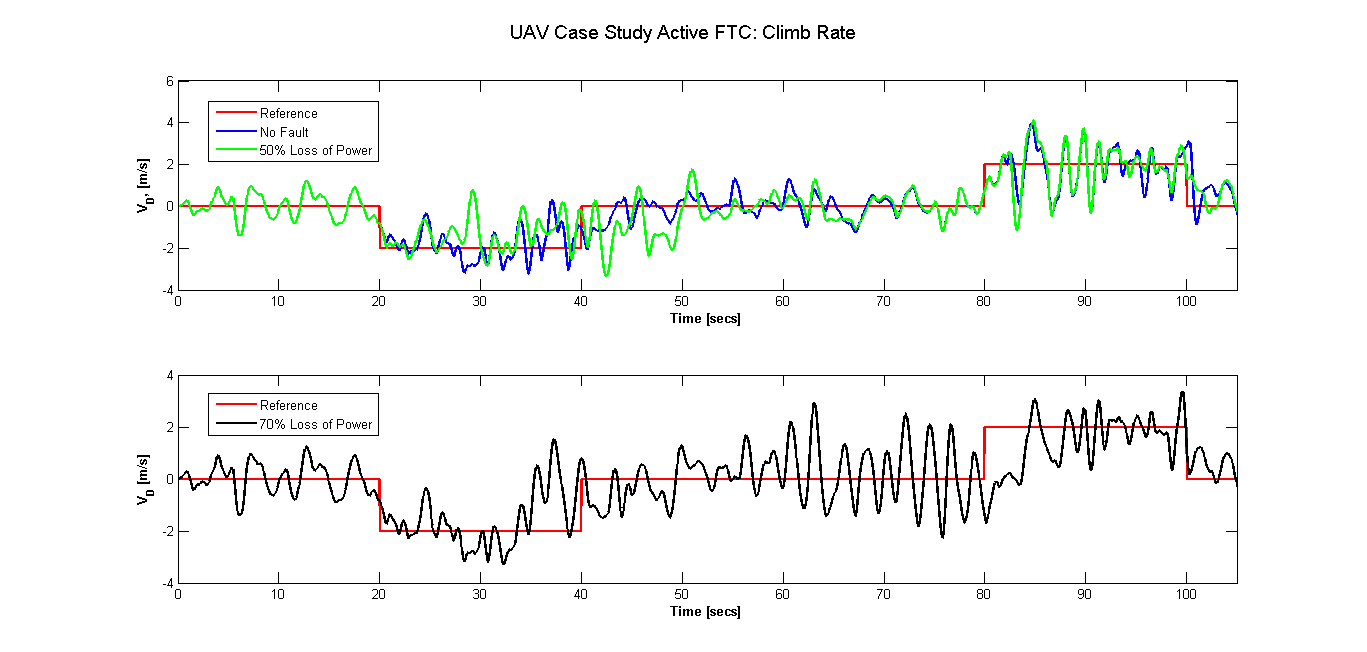} 
%\vspace{-0.5in}
\caption{UAV Case Study Active FTC: Climb Rate, $V_D$}
\label{fig:chap6_activeFTC_controls_Vd_subs}
\end{figure}

The climb rate (or vertical speed) response is shown in figures \ref{fig:chap6_activeFTC_controls_Vd} and \ref{fig:chap6_activeFTC_controls_Vd_subs} and shows that during power loss the vertical speed oscillates between the upper and lower constraints values.  This is to be expected as the elevator is working harder to regulate the speed, and hence bounces between the two limits. 

\begin{figure}[H]
%\hspace{-1.2in}
\center
\includegraphics[scale=0.35]{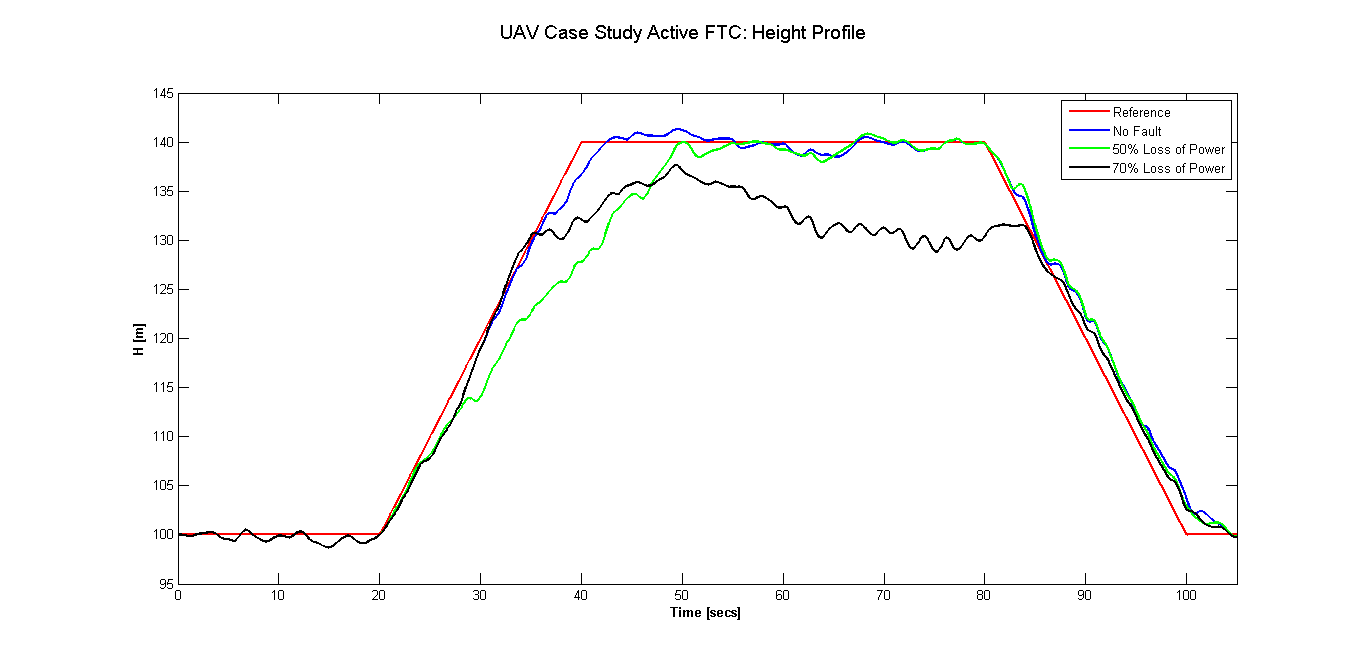} 
%\vspace{-0.5in}
\caption{UAV Case Study Active FTC: Height Profile}
\label{fig:chap6_activeFTC_controls_Height}
\end{figure}

Finally the trajectories flown by the aircraft in each scenario are shown in figure \ref{fig:chap6_activeFTC_controls_Height}.  The interesting points to note here are that in the highest loss of power case (scenario 3) the aircraft did not have enough power to reach the highest flight altitude so instead cruised at an altitude it was capable of flying.  Once the straight and level phase of the flight was over the controller was able to fly the aircraft back onto the demanded trajectory.  The $50\%$ loss of power case shows that aircraft still had enough power to fly back onto the demanded path halfway through the straight and level at altitude flight phase.

\section{Conclusion}\label{section:chap6_sof}
In this paper we demonstrated the successful application of an NMPC based active FTC design on an actual UAV which is currently in operation.  The controller was implemented using actual aircraft data and shows great promise for active fault tolerant flight control.  

\section*{Appendix}

An Appendix, if needed, appears before the acknowledgements.

\section*{Acknowledgements}
An Acknowledgements section, if used, immediately precedes the References. Sponsorship and financial support acknowledgements should be included here.

\end{document}